\documentclass[12pt]{article}
\usepackage{amssymb}
\usepackage[leqno]{amsmath}

\usepackage{amsfonts}
\usepackage{amstext}
\usepackage{amsthm}
\usepackage{bm} 

\newcommand{\Aaa}{{\cal{A}}}

\newcommand{\Ef}{{\cal{F}}}
\newcommand{\Gee}{{\cal{G}}}
\newcommand{\Haa}{{\cal{H}}}

\newcommand{\Be}{{\Bbb{B}}}

\newtheorem{definition}{{Definition}}[section]
\newtheorem{theorem}[definition]{{Theorem}}
\newtheorem{corollary}[definition]{{Corollary}}

\newtheorem{lemma}[definition]{\noindent {Lemma}}

\newtheorem{claim}[definition]{\noindent {Claim}}
\theoremstyle{definition} 

\newtheorem{remark}[definition]{\noindent{Remark}}

\newcommand{\inv}[2]{{#1}^{-1}[{#2}]}

\def\eqdef{\mbox{\bf\ :=\ }} 

\newcommand{\rfs}[1]{{\rm #1}}
\newcommand{\fnn}[3]{#1:#2 \rightarrow #3}

\newcommand{\setn}[2]{\{#1\: :\:#2\}}
\newcommand{\setof}[2]{\{#1\: :\:#2\}}

\newcommand{\seqn}[2]{\langle#1\: :\:#2\rangle }

\newcommand{\sngltn}[1]{ \{#1\} }

\newcommand{\pair}[2]{\langle #1 , \:#2 \rangle }
\newcommand{\triple}[3]{\langle\kern1pt#1 , \:#2 , 
      \:#3 \kern1pt\rangle }
\newcommand{\trpl}[3]{\langle\kern1pt#1 , \:#2 , 
      \:#3 \kern1pt\rangle }
\newcommand{\forthuple}[4]{\langle\kern1pt#1 , \:#2  , 
      \:#3  , \:#4\kern1pt\rangle }
\newcommand{\qdrpl}[4]{\langle\kern1pt#1, #2, #3, 
   #4\kern1pt\rangle}
\newcommand{\fifthtpl}[5]{\langle\kern1pt#1, #2, #3, 
   #4, #5\kern1pt\rangle}
\newcommand{\sixtpl}[6]{\langle\kern1pt#1, #2, #3, 
   #4, #5, #6 \kern1pt\rangle}

\newcommand{\meet}{\wedge}
\newcommand{\join}{\vee}

\renewcommand{\Bbb}{\mathbb}

\newenvironment{pf}{\begin{proof}}{\end{proof}}

\providecommand{\cal}{\mathcal}

\def\calW{\mathcal{W}}

\newcommand{\sett}[2]{{\{#1\}}_{#2}}
\newcommand{\sn}[1]{\{#1\}} 
\newcommand{\loe}{\le}
\newcommand{\goe}{\ge}

\newcommand{\vsum}{\vec{\sum}}
\newcommand{\vdot}{\;\vec{\cdot}\;}

\newcommand{\congtop}{\cong^{\rfs{top}}} 

\newcommand{\cbd}[2]{{\rfs D}^{#1}({#2})} 
\newcommand{\rk}{\operatorname{rk}} 
\newcommand{\map}[3]{#1\colon #2 \to #3} 
\newcommand{\tplank}{{\mathbb T\mathbb P}} 
\newcommand{\gplank}{{\mathbb G\mathbb P}} 
\renewcommand{\coprod}{\uplus} 

\newcommand{\eps}{\varepsilon}
\newcommand{\fs}[1]{{\rfs{FS}\left(#1\right)}} 
\newcommand{\F}[1]{{\widehat{F}\left(#1\right)}}

\newcommand{\img}[2]{#1\left[#2\right]} 
\newcommand{\power}[1]{{\mathcal P}(#1)}
\newcommand{\ind}{\rfs{Ind}} 
\providecommand{\cf}{\operatorname{cf}}
\newcommand{\calw}{{\calW^*_2}} 
\newcommand{\seq}[1]{\langle #1 \rangle}
\begin{document}
\thispagestyle{empty}
\phantom{bfk}
\vspace{-20mm}
\begin{center}
{\Large\bf Free Boolean algebras over unions 

of two well orderings\vspace{3mm}}

{\bf Robert Bonnet}%
\footnote{This work is partially supported by 
the Institute of Mathematics, Akademia \'Swi\c{e}tokrzyska, Kielce,
August - September 2006.}

{\small Laboratoire de Math\'ematiques,

Universit\'e de Savoie, Le Bourget-du-Lac, France\vspace{1mm}
}

{\bf Latifa Faouzi}%
\footnote{This work is supported by 
the Institut de Math\'{e}matique de Luminy,
(CNRS-UMR~6206) and by the 
Institut Camille Jordan, Lyon (CNRS-UMR~5208), September - October 2005.}

{\small 
\smallskip
{\small U.F.R. Math\'{e}matiques Discr\`{e}tes et Applications, F.S.T.}
\\
Universit\'{e} Sidi Mohamed Ben Abdellah,
F\`es, Maroc\vspace{1mm}
}

{\bf Wies{\l}aw Kubi\'s}%
\footnote{This work is partially supported by 
the Laboratoire de Math\'{e}matiques (CNRS-UMR 5127)
as invited Professor,
Universit\'{e} de Savoie,  July~2006, and 
the Fields Institute for Research in Mathematical Sciences, 
Toronto, January - March 2007.
}

{\small 
Instytut Matematyki, Akademia \'Swi\c{e}tokrzyska,

Kielce, Poland
}
\end{center}
\vspace{-4mm}
\begin{abstract}
\noindent
Given a partially ordered set $P$ there exists the most general Boolean algebra $\F P$
which contains $P$ as a generating set, called the {\it free Boolean algebra} over $P$. 
We study free Boolean algebras over posets of the form $P=P_0\cup P_1$, 
where $P_0,P_1$ are well orderings. 
We call them {\it nearly ordinal algebras}.

Answering a question of Maurice Pouzet, we show that 
for every uncountable cardinal $\kappa$ there are $2^\kappa$ 
pairwise non-isomorphic nearly ordinal algebras of cardinality $\kappa$. 

Topologically, free Boolean algebras over posets correspond to compact 
$0$-dimensional distributive lattices. In this context, we classify all closed sublattices 
of the product $(\omega_1+1)\times(\omega_1+1)$, thus showing that 
there are only $\aleph_1$ many of them.
In contrast with the last result, we show that there are $2^{\aleph_1}$
topological types of closed subsets of the 
Tikhonov plank $(\omega_1+1)\times(\omega+1)$. 
\end{abstract}
\begin{footnotesize}
{\bf Keywords:} Well quasi orderings, Poset algebras, 
Superatomic Boolean algebras,
Compact distributive lattices.
\newline
{\bf Mathematics Subject  Classification 2000 (MSC2000)}
\newline
{\bf Primary:}
03G05,  
06A06.  
{\bf Secondary:}
06E05,  
08A05,  
54G12.  
\medskip
\newline
{\bf E-mail:} 
bonnet@in2p3.fr (R. Bonnet),
latifaouzi@menara.ma (L. Faouzi),
\newline%
\phantom{\bf E-mail:} 
\,wkubis@pu.kielce.pl (W. Kubi\'s).
\end{footnotesize}
\newpage
\tableofcontents
\setcounter{footnote}{0}
\newpage
\section{Introduction}
\label{section-intro}
Assume that a given Boolean algebra $\Be$ is generated by a well ordered subset $P$.
Then $\Be$ is an interval algebra whose Stone space is homeomorphic to $\alpha+1$,
where $\alpha$ is the order type of $P$. We call $\Be$ an {\it ordinal algebra} and we
call its Stone space an {\it ordinal space}. It is clear that for a given infinite cardinal 
$\kappa$, there are only $\kappa^+$ many types of ordinal algebras 
of cardinality $\kappa$. 

Given a poset $P$, we denote by $\F P$ the {\it free Boolean algebra} over $P$ 
(the precise definitions are given in the next section). 
We shall call $\F P$ a {\it nearly ordinal algebra} if $P=P_0\cup P_1$, 
where $P_0$, $P_1$ are well ordered sets.
We show that for every uncountable cardinal $\kappa$ there exist $2^{\kappa}$
isomorphic types of nearly ordinal algebras of cardinality $\kappa$. 
This answers a question of Maurice Pouzet \cite{pouzet}, 
originally asked in a more general context of posets 
which are well founded and narrow (the so called {\it well quasi-ordered} posets).

One of the simplest Boolean algebras which are not generated by a chain is 
$\Be = \F{\omega_1\coprod \omega}$, where $P_0\coprod P_1$ 
denotes the disjoint union of posets $P_0$, $P_1$ with no extra relations.
Topologically, $\Be$ corresponds to the {\it Tikhonov plank} 
$\tplank = (\omega_1+1)\times(\omega+1)$. 
We show that $\tplank$ has $2^{\aleph_1}$ many isomorphic types of closed sets. 
On the other hand, we show that $\tplank$ has only countably many topological types 
of uncountable closed sets which are at the same time sublattices of $\tplank$. 
Finally, we classify topological types of compact sublattices of $(\omega_1+1)^2$. 
We show that an uncountable closed sublattice of $(\omega_1+1)^2$ is homeomorphic
to one of the following spaces: $(\omega_1+1)\times(\alpha+1)$, where 
$\alpha\loe\omega_1$, and $\triangle = \setof{\pair xy\in(\omega_1+1)^2}{x\loe y}$. 
This implies that there are only $\aleph_1$ isomorphic types of Boolean algebras 
of the form $\F{P_0\cup P_1}$, where $P_0,P_1$ are well ordered chains of order type 
$\loe\omega_1$.
\section{Terminology and basic facts}%
\label{section-teminology}%
All topological spaces considered in this paper are Hausdorff and zero-dimen\-sio\-nal.

Given a compact space $X$, we denote by $\cbd \alpha X$ the $\alpha$-th 
Can\-tor-Ben\-dix\-son derivative of $X$. 
Specifically, $\cbd 1X$ is the set of all non-isolated points of $X$, 
$\cbd{\alpha+1}X=\cbd 1{\cbd\alpha X}$ and 
$\cbd \delta X = \bigcap_{\xi<\delta}\cbd \xi X$ for any limit $\delta$.
If $X$ is scattered, then the {\it Cantor-Bendixson rank} of $X$, denoted by $\rk(X)$, 
is the first ordinal $\alpha$ such that $\cbd{\alpha+1} X$ is empty.
A compact scattered space $X$ is {\it unitary} whenever $\cbd{\rk(X)}X$
has a unique element, called the {\it end-point} of $X$ and denoted by $e(X)$.
For $x \in X$, we denote by $\rfs{rk}^X(x)$ the unique ordinal
$\theta$ such that $x \in \cbd \theta X \setminus \cbd {\theta+1} X$
and $\rfs{rk}^X(x)$ is called the {\it Cantor-Bendixson rank of $x$ in $X$}.

Let  $X$ and $Y$ be topological spaces.
We denote by $X \congtop Y$ the fact that $X$ and $Y$ are homeomorphic.

A compact space $X$ is {\it retractive} if every closed subset of $X$ is a retract of $X$,
i.e. for every closed $Y \subseteq X$ there is a continuous map $\map fXY$ such that 
$f(x)=x$ for every $x\in Y$.
Let us mention a theorem of M. Rubin~\cite[Theorem~5.1]{R} 
(see also \cite[\S15, Theorem 15.18]{K}):
every Boolean space which is a continuous image of a linearly ordered compact is 
retractive.

The {\it Tikhonov plank} is the space $\tplank \eqdef (\omega_1+1)\times(\omega+1)$,
where both factors are considered as linearly ordered compacta. 
It is straight to see that $\tplank$ is not retractive, namely the closed set 
$$A=\Bigl((\omega_1+1)\times\sn\omega\Bigr) \cup \Bigl(\sn{\omega_1}\times(\omega+1)\Bigr)$$
is not a retract of $\tplank$. 
In particular, $\tplank$ is not a continuous image of any linearly ordered compact space.

Given a partially ordered set $P$, we denote by $\F P$ 
the {\it free Boolean algebra over} $P$. 
More precisely, $\F P$ is a Boolean algebra together with 
an order preserving embedding $\map {i_P}P{\F P}$ 
such that for every order preserving map $\map{f}{P}{\Be}$ into an arbitrary 
Boolean algebra $\Be$ there exists a unique homomorphism $\map{h}{{\F P}}{\Be}$ 
such that $h \circ i_P=f$. 
These two properties define $\F P$ uniquely. 
Topologically, $\F P$ can be defined to be the clopen algebra of the space 
$$\fs P \eqdef \setof{F\subseteq P}
{(\forall\;x,y\in P)\;
(x\in F \mbox{\rm\ and } x\loe y \mbox{\rm\ imply } y\in F\,)},$$
with the topology inherited from the Cantor cube $\power P$ (the powerset of $P$).
The embedding $i_P$ is defined by the formula $i_P(p)=\setof{F\in\fs P}{p\in F}$.
Elements of $\fs P$ are called {\it final segments} of $P$. 
Note that $\fs P$ has a natural structure of a lattice of sets: 
the union as well as the intersection of any family of final segments is again a 
final segment. 
It turns out that the original poset can always be reconstructed from the lattice structure
as the set of all clopen prime filters with inclusion. 
On the other hand, every compact $0$-dimensional distributive lattice 
$\seq {L,0,1,\join,\meet,\tau}$ is isomorphic 
(in the category of topological lattices) 
to $\seq{\fs P, \emptyset, P, \cup,\cap, \tau_P}$ for some poset $P$, 
where $\tau_P$ is the topology of $\fs P$.
In fact, $P\mapsto \fs P$ is a contravariant functor from the category of posets 
onto the category of compact $0$-dimensional distributive lattices, 
whose inverse is defined by assigning to a lattice the family of all its 
clopen prime filters. 
We refer to \cite[\S2]{ABKR} for more details.

Given a compact lattice $L$, we shall denote by $1_L$ and $0_L$ 
the maximal and the minimal element of $L$ respectively.
For our purposes, it is more convenient to consider the lattice $\fs P$ 
with reversed inclusion, i.e. the join of $x,y\in\fs P$ is $x\cap y$ and 
their meet is $x\cup y$. 
In particular, $1_{\fs P}$ is the empty final segment 
and $0_{\fs P}$ is the whole poset $P\in\fs P$. 
For example, if $P$ is a well ordered chain then $\pair{\fs P}{\supseteq}$ 
is a well ordered compact chain. 
Note that a basic neighborhood of a point $x\in \fs P$ is of the form 
$s^+\cap t^-$, where $s,t\subseteq P$ are finite sets, 
$s^+  \eqdef \setof{y\in \fs P}{s\subseteq y}$ 
and $t^-\eqdef \setof{y\in \fs P}{t\cap y=\emptyset}$.
Observe that $0_P=P$ is an isolated point of $\fs P$, 
whenever there is a finite set $s\subseteq P$ such that the final segment 
generated by $s$ is $P$. 
In that case $\sn P=s^+$. 
This happens, for example, in case where $P$ is the union of finitely 
many well ordered chains.

Given posets $P_0$, $P_1$ we denote by $P_0\coprod P_1$ the disjoint (incomparable)
sum of $P_0$, $P_1$, i.e. $P_0\coprod P_1$ is the disjoint union of $P_0$, $P_1$ 
and $x\loe y$ in $P_0\coprod P_1$ iff $x,y\in P_i$ and $x\loe y$ in $P_i$ for some $i<2$.
We denote by $P_0\cup P_1$ a poset wich is covered by
two posets $P_0$ and $P_1$.
Recall that a {\it nearly ordinal algebra} is a Boolean algebra of the form 
$\F{P_0\cup P_1}$, where $P_0,P_1$ are well ordered chains. 
We shall call $\fs{P_0\cup P_1}$ a {\it nearly ordinal space}.
\section{The number of nearly ordinal Boolean algebras}
In this section we prove the announced result on the number of nearly ordinal 
Boolean algebras of a given cardinality. 
Before we state the precise statement, we shall introduce some necessary notions.

Given unitary scattered compact spaces $X$ and $Y$, we say that $X$ and $Y$ 
are {\it almost homeomorphic} if there are clopen neighborhoods $U$, $V$ of 
the end-points of $X$ and $Y$ respectively such that $U$ and $V$ are homeomorphic. 

We now define for each uncountable cardinal $\kappa$ an ordinal $r(\kappa)$ as follows.
If $\kappa>\aleph_1$ is a successor cardinal then 
we set $r(\kappa) = \kappa$. 
If $\kappa$ is a limit cardinal which is not a strong limit (i.e. $2^\lambda>\kappa$ 
for some $\lambda<\kappa$) then again we set $r(\kappa)=\kappa$. 
If $\kappa$ is a strong limit 
(i.e. $\kappa$ is a limit ordinal and $2^{\lambda}<\kappa$ for $\lambda<\kappa$) 
then we set $r(\kappa)=\kappa+\kappa$. It remains to define $r(\aleph_1)$. 
Let
$$E=\setof{\theta<\omega_2}{\omega^\theta=\theta\text{ and }\cf\theta=\omega_1}.$$
Then $|E|=\aleph_2$. 
Let $E_0$ be the initial segment of $E$ whose order type is $\omega_1$. 
Finally, let $\varpi_1 = \sup(E_0)$ and set $r(\aleph_1)=\varpi_1+\omega_1$.

Note that always $\kappa\loe r(\kappa)<\kappa^+$ and in fact $\kappa<r(\kappa)$ only if $\kappa=\aleph_1$ or $\kappa$ is a strong limit cardinal.
\begin{theorem}
\label{thm-3.1}
Let $\kappa$ be an uncountable cardinal. 
There exists a family $\seqn{X_\alpha}{\alpha<2^\kappa}$ of nearly ordinal spaces with 
the following properties:
\begin{enumerate}
	\item[\rm(1)] Each $X_\alpha$ is unitary of cardinality $\kappa$ and 
	of Cantor-Bendixson rank $r(\kappa)$. 
	\item[\rm(2)] For each $\alpha$, the space $X_\alpha$ is not retractive. Specifically: every neighborhood of $e(X_\alpha)$ contains a copy of the Tikho\-nov plank.
	\item[\rm(3)] If $\alpha\ne \beta$ then $X_\alpha$ and $X_\beta$ are not 
	almost homeomorphic.
\end{enumerate}
In particular, for each uncountable cardinal $\kappa$ there exist $2^\kappa$ pairwise non-isomorphic nearly ordinal Boolean algebras of size $\kappa$.
\end{theorem}
The rest of this section is devoted to the proof of Theorem~\ref{thm-3.1}. We start with some notation and definitions.

We shall denote by $\ind$ the class of all indecomposable ordinals. 
Recall that an ordinal $\rho>0$ is {\it indecomposable} if $\alpha+\beta<\rho$ 
whenever $\alpha,\beta<\rho$.
In this case the linearly ordered compact space $\rho+1$ is unitary and $e(\rho+1)=\rho$. 
The Cantor-Bendixson rank of $\rho+1$ will be denoted by ${\ln}(\rho)$.

We shall denote by $\calw$ the class of all lattices isomorphic to $\fs P$, 
where $P$ is covered by two well ordered chains and $\emptyset=1_\fs P$ 
is the unique element of maximal Cantor-Bendixson rank in $\fs P$. 
So every lattice in $\calw$ is a unitary scattered compact space.

Let $\rho>0$ be an ordinal and let $\seqn{L_\xi}{\xi<\rho}$ be a collection of 
compact lattices, where $L_\xi = \fs{P_\xi}$ for some poset $P_\xi$. 
Let $P=\sum_{\xi<\rho}P_\xi$ by the lexicographic sum of $\setof{P_\xi}{\xi<\rho}$. 
We shall write $\vsum_{\xi<\rho}L_\xi$ for the lattice $\fs{P}$. 
In case where $L_\xi=L$ for every $\xi$, we shall write $L\vdot\rho$ 
instead of $\vsum_{\xi<\rho}L_\xi$.

Let $K=\vsum_{\xi<\rho}L_\xi$, where each $L_\xi$ is a well founded distributive lattice
isomorphic to $\fs{P_\xi}$, where $P_\xi$ is a union of two well ordered chains. 
Observe that $K$ is a nearly ordinal space, because $\sum_{\xi<\rho}P_\xi$ is again 
a union of two well ordered chains. 
Recall that we consider $\fs P$ with the reversed ordering. 
Observe that, given $\eta<\rho$ we can identify $x\in\fs{P_\eta}$ with 
$$\widehat{x} \eqdef \setof{p\in\sum_{\xi<\rho}P_\xi}{(\exists\;q\in x)\;(p \goe q)}
\; \in \; \fs P.$$
By this way we get a natural embedding of $L_\eta$ into $K$. 
Notice that the empty final segment of $P_\eta$ is actually identified with 
the full final segment of $P_{\eta+1}$. 
In other words, $1_{L_\eta}$ is identified with $0_{\eta+1}$. 
Further, observe that $L_\eta\setminus\sn{0_{L_\eta}}$ is clopen in $K$. 
Indeed, denoting by $s_\alpha$ the set of minimal elements of $P_\alpha$, we have
$$  L_ \eta^- \eqdef L_\eta\setminus\sn{0_{L_\eta}} 
= \setof{x\in \fs P}{s_\eta\not\subseteq x\text{ and }s_{\eta+1}\subseteq x}.$$
The set on the righ-hand side is evidently clopen, because each $s_\alpha$ is a finite set.
Finally, observe that for a limit ordinal $\delta<\rho$, 
a basic neighborhood of $0_{L_\delta}$ is of the form 
$V(\alpha,\delta)
\eqdef
\sn{0_{L_\delta}} \cup L_\alpha^- \cup \bigcup_{\xi\in(\alpha,\delta)}L_\xi$
%
where $\alpha<\delta$. 
Indeed, $V(\alpha,\delta) = s_\delta^+ \setminus s_\alpha^+$ and thus $V(\alpha,\delta)$  is clopen in $K$; 
and a typical neighborhood of $0_{L_\delta}$  is of the form $s^+\cap t^-$, 
where $s,t$ are finite subsets of $P$. 
Necessarily $s\subseteq \bigcup_{\xi\goe\delta}P_\xi$. 
Thus $0_{L_\delta} \in V(\alpha,\delta)\subseteq s^+\cap t^-$, 
where $\alpha$ is such that $t\subseteq\bigcup_{\xi<\alpha}P_\xi$. 
Similarly, sets of the form 
$V(\alpha,\rho) = \sn{1_L}  \cup L_\alpha^- \cup \bigcup_{\xi>\alpha}L_\xi$, 
where $\alpha<\rho$, form an open base at $1_L=\emptyset$.
\begin{lemma}
\label{lemma-3.2}
Let $\rho,\theta>0$ be ordinals and assume $\rho$ is indecomposable. 
Further, let $\setof{Y_\alpha}{\alpha<\rho}\subseteq \calw$ be such that 
$\rk(Y_\alpha) = \theta$ for each $\alpha<\rho$. 
Let $Y=\vsum_{\alpha<\rho}Y_\alpha$. 
Then 
$Y\in\calw$, $\cbd\theta Y\setminus \cbd{\theta+1}Y=\setof{1_{Y_\xi}}{\xi<\rho}$ 
and the Cantor-Bendixson rank of $Y$ equals $\theta+\ln\rho$.
\end{lemma}
\begin{pf}
Recall that $Y_\xi\setminus\sn{0_{Y_\xi}}$ is a clopen subset of $Y$ 
and $e(Y_\xi)=1_{Y_\xi}$, because the minimal element of a well founded distributive
lattice is isolated. 
Thus the Cantor-Bendixson rank of $e(Y_\xi)=1_{Y_\xi}$ in $Y$ is the same as 
its rank in $Y_\xi$. 
It follows that $\setof{1_{Y_\xi}}{\xi<\rho}\subseteq\cbd\theta{Y}$. 
By the same reason, $\cbd\theta Y\cap Y_\xi \subseteq \{0_{Y_\xi}, 1_{Y_\xi}\}$. 
Thus
$$\cbd \theta Y = \setof{1_{Y_\xi}}{\xi<\rho} \cup 
\setof{ 0_{Y_\eta} }{ \eta\text{ is a limit ordinal }<\rho} \cup \sn{1_Y} \, ,$$
because each $1_{Y_\xi}$ is identified with $0_{Y_{\xi+1}}$. 
Hence $\cbd\theta Y$ is homeomorphic to $\rho+1$. 
Finally, $\rk(\rho+1)=\ln\rho$ and $e(Y)=1_Y$, which shows that $Y\in\calw$.
\end{pf}
\begin{lemma}
\label{lemma-3.3}
Let $\rho,\theta>0$ be ordinals and assume $\rho$ is indecomposable. 
Further, let $\setof{Y_\alpha}{\alpha<\rho}\subseteq\calw$ and 
$\setof{Z_\alpha}{\alpha<\rho}\subseteq\calw$ 
be such that $\rk(Y_\alpha)=\rk(Z_\alpha)=\theta$ for each $\alpha<\rho$. 
Let $Y=\vsum_{\alpha<\rho}Y_\alpha$\,, \,$Z=\vsum_{\alpha<\rho}Z_\alpha$ and 
assume that $Y$ and $Z$ are almost homeomorphic. 
Then there exists $\alpha<\rho$ such that for every $\beta\goe\alpha$ 
there exists $\xi$ with the property that $Y_\beta$ is almost homeomorphic to $Z_\xi$.
\end{lemma}
\begin{pf}
Let $\map hUV$ be a homeomorphism between clopen neighborhoods of $1_Y$ 
and $1_Z$ respectively.
Let $\alpha<\rho$ be such that $Y_\beta\subseteq U$ for every $\beta\goe\alpha$. 
Fix $\beta\goe\alpha$. 
By assumption and by Lemma~\ref{lemma-3.2}, we know that $1_{Y_\xi}$ are 
precisely the elements of $Y$ whose Cantor-Bendixson rank is $\theta$. 
The same applies to~$Z$. 
Thus $h(1_{Y_\beta}) = 1_{Z_\xi}\in V$ for some $\xi$. 
Let $U'=U\cap(Y_\beta\setminus\sn{0_{Y_\beta}}$. 
Then $U'$ is a clopen neighborhood of $1_{Y_\beta}$ in $Y_\beta$.
Let $V'=\img h{U'}\cap(Z_\xi\setminus\sn{0_{Z_\xi}}$. 
Then $V'$ is a clopen neighborhood of $1_{Z_\xi}$ in $Z_\xi$, 
homeomorphic to $\inv h{V'}\subseteq U'\subseteq Y_\beta$. 
Thus $Y_\beta$ is almost homeomorphic to $Z_\xi$, which completes the proof.
\end{pf}
\begin{lemma}
\label{lemma-3.4}
Let $\kappa$ be an uncountable cardinal and let $\gamma$ be an ordinal such that $\gamma+\kappa = r(\kappa)$.
Assume further that $\setof{ Y^i }{ i<\kappa }$ is a family of spaces satisfying the following conditions:
\begin{enumerate}
	\item[$\star_1(\kappa)$] $Y^i$ is a member of $\calw$.
	\item[$\star_2(\kappa)$] $| Y^i |  \loe \kappa$.
	\item[$\star_3(\kappa)$] $\rk(Y^i) = \gamma$.
	\item[$\star_4(\kappa)$] $Y^i$ is unitary.
	\item[$\star_5(\kappa)$] $Y^i$ and $Y^j$ are not almost homeomorphic for distinct $i , j < \kappa$.
	\item[$\star_6(\kappa)$] If $U$ is a neighborhood of $e(Y^i)$ in $Y^i$ then $U$ contains a closed subspace homeomorphic to the Tikhonov plank.
\end{enumerate}
Then there exists a family of spaces $\setof{X_\alpha}{\alpha<2^\kappa}\subseteq\calw$
satisfying the assertions of Theorem~\ref{thm-3.1}.
\end{lemma}
\begin{pf}
We construct the family $\setof{ X_{\alpha} }{ \alpha<2^\kappa }$ as follows.
Let $\setof{ A_\alpha }{ \alpha<2^\kappa }$ be an enumeration of
all subsets $A \subseteq \kappa$ such that $| A |Ê\geq 2$.
For each $\alpha<2^\kappa$, let $\fnn{f_\alpha}{\kappa}{ A_\alpha }$
be such that $\inv{f_\alpha}{\xi}$ is cofinal in $\kappa$ for every $\xi \in A_{\alpha}$.
Let 
$X_{\alpha} \eqdef \vsum_{ \zeta < \kappa } Y^{f_\alpha(\zeta)} $.
We prove that the family $\setof{ X_{\alpha} }{ \alpha<2^\kappa }$ 
is as required.

It is clear that $X_\alpha$ is a unitary scattered compact of cardinality $\kappa$. 
By Lemma~\ref{lemma-3.2}, the Cantor-Bendixson rank of $X_\alpha$ 
equals $\gamma+\kappa=r(\kappa)$.
Now let $U$ be a clopen neighborhood of $e(X_{\alpha})$ in $X_{\alpha}$.
There is $\zeta<\kappa$ such that $Y^{f_\alpha(\zeta)} \subseteq X_{\alpha}$.
By $\star_6(\kappa)$, $Y^{f_\alpha(\zeta)}$ contains a copy of~$\tplank$.
This shows condition (2) of Theorem~\ref{thm-3.1}. It remains to prove~(3).

Fix $\alpha\ne\beta$. We need to show that $X_\alpha$ and $X_\beta$ 
are not almost homeomorphic.
Assume for instance that $\xi \in A_{\beta} \setminus A_{\alpha}$.
We show that the assertion of Lemma~\ref{lemma-3.3} is not fulfilled.
For fix any $\delta<\kappa$.
Since $\inv{f_\beta}\xi$ is cofinal in $\kappa$, 
let $\zeta > \delta$ be such that $f_\beta(\zeta) = \xi$.
On the other hand, $\xi \not\in A_{\alpha}$, which means that  
\,$f_\alpha(\eta) \neq \xi$ for every $\eta<\kappa$.
Hence, by~$\star_5(\kappa)$, 
$Y^{f_\beta(\zeta)}$ is not homeomorphic to $Y^{f_\alpha(\eta)}$ for every $\eta<\kappa$.
Since $\zeta>\delta$, by Lemma~\ref{lemma-3.3}, we conclude that 
$X_{\alpha}$ and $X_{\beta}$ are not almost homeomorphic.
\end{pf}
Theorem~\ref{thm-3.1} is a consequence of the following statement combined with Lemma~\ref{lemma-3.4}.
\begin{lemma}
\label{lemma-3.5}
For every uncountable cardinal $\kappa$ there exist a family $\seqn{ Y^i }{ i<\kappa }$ 
and an ordinal $\gamma$ with $\gamma+\kappa=r(\kappa)$, 
satisfying conditions $\star_1(\kappa)$ -- $\star_6(\kappa)$ of Lemma~\ref{lemma-3.4}.
\end{lemma}
\begin{pf} We use induction on the cardinal $\kappa$. 
Assume the statement has been proved for all uncountable cardinals $\lambda<\kappa$.
We consider three cases.
\medskip
\newline
{\bf Case 0.}
$\kappa=\aleph_1$.
\medskip
\newline
We set $\gamma=\varpi_1$. 
Given a limit ordinal $\theta$, let $K_\theta=(\theta+1)^2$ and let $T_\theta$ be 
the space obtained from two copies of $\theta+1$ by identifying their last elements. 
In other words, $T_\theta$ can be described as the linearly ordered space 
$\theta+1+\theta^*$, where $\theta^*$ is the set $\theta$ with the reversed ordering. 
The following claim is trivial, after noting that $\rk(\theta+1)=\theta$ 
and $e(\theta+1)=\theta$, for any $\eps$-ordinal $\theta$.
\begin{claim}
\label{claim-3.6}
Let $\theta>0$ be an ordinal such that $\omega^\theta=\theta$. 
Then
\begin{enumerate}
	\item[\rm(a)] $\rk(T_\theta) = \theta$ and $T_\theta$ is unitary.
	\item[\rm(b)] $T_\theta$ is not an ordinal space.
	\item[\rm(c)] $\cbd{(\theta)}{K_\theta} \congtop T_\theta$.
	\item[\rm(d)] $\rk(K_\theta) = \theta \cdot 2$ and
	$e(K_ \theta) = e(T_\theta) = 1_{K_\theta}$.
\end{enumerate}
\end{claim}
Recall that we have denoted by $E_0$ the set consisting of the first $\omega_1$ 
many $\eps$-ordinals of uncountable cofinality. 
In particular, $\omega^\theta=\theta$ for $\theta\in E_0$. 
Given $\theta\in E_0$, define $X_\theta = K_\theta \vdot \varpi_1$. 
Finally, fix a one-to-one enumeration $\seqn{\theta_i}{i<\omega_1}$ of $E_0$ 
and define $Y^i=X_{\theta_i}$. 

Since $\fs {\theta\uplus \theta} = K_\theta$, 
clearly, $X_\theta\in\calw$ is unitary and of cardinality $\aleph_1$. 
Further, every neighborhood of $e(X_\theta) = 1_{X_\theta}$ contains 
some $K_\theta$ and $\tplank$ is a sublattice of $K_\theta$, 
because $\omega_1\loe\theta$. 
This shows that the family $\seqn{Y^i}{i<\omega_1}$ satisfies $\star_i(\aleph_1)$ 
for $i=1,2,4,6$.
By Claim~\ref{claim-3.6}(d) together with Lemma~\ref{lemma-3.2}, 
we have that 
$$\rk(X_\theta)=\theta\cdot 2 + \ln\varpi_1 = \varpi_1=\gamma.$$
This shows $\star_3(\aleph_1)$. 
It remains to show $\star_5(\aleph_1)$.
For fix $\theta<\theta'$ in $E_0$. 
Then $\theta\cdot2 <\theta'\cdot 2$, so by Claim~\ref{claim-3.6}(d), 
$K_\theta$ and $K_{\theta'}$ cannot be almost homeomorphic, 
because of the Cantor-Bendixson rank. 
Hence $X_\theta$ and $X_{\theta'}$ are not almost homeomorphic by 
Lemma~\ref{lemma-3.3}.
Thus $\seqn{Y_i}{i<\omega_1}$ satisfies $\star_1(\omega_1)$ -- $\star_6(\omega_1)$.
\medskip
\newline
{\bf Case 1.}
$\kappa$ satisfies:
$\lambda<\kappa$ and $2^\lambda\goe\kappa$ for some $\lambda>\aleph_0$.
\medskip
\newline
So, either $\kappa$ is a successor cardinal or else $\kappa$ is a limit, 
but not a strong limit cardinal.
By the induction hypothesis applied to $\lambda$, followed by Lemma~\ref{lemma-3.4},
there is a sequence $\setof{ Y^i }{ i < 2^\lambda}$ satisfying the assertions of 
Theorem~\ref{thm-3.1}.

Note that $| Y^i | = \lambda$ and $\rk(Y^i) = r(\lambda)$.
Let $\gamma=r(\lambda)<\lambda^+\loe \kappa$.
It is clear that the family $\setof{ Y^i }{ i<\kappa }$ satisfies conditions 
$\star_1(\kappa)$ -- $\star_6(\kappa)$.
\medskip
\newline
{\bf Case 2.}
$\kappa$ is a strong limit cardinal.
\medskip
\newline
Fix a strictly increasing sequence of infinite cardinals 
$\seqn{\kappa_\alpha}{\alpha<\cf\kappa}$ 
with $\kappa= \sup_{\alpha<\cf\kappa}\kappa_\alpha$. 
For each $\alpha<\cf\kappa$ let $\Ef_\alpha\subseteq\calw$ be a family satisfying 
$\star_1(\kappa_\alpha^+)$ -- $\star_6(\kappa_\alpha^+)$, 
obtained by the induction hypothesis. 
Let $\rho_\alpha$ denote the common Cantor-Bendixson rank of the spaces 
from $\Ef_\alpha$. 
We define families $\Gee_\alpha \subseteq \Ef_\alpha$ satisfying 
the following conditions.
\begin{enumerate}
	\item[(i)] $|\Gee_\alpha|=\kappa_\alpha^+$.
	\item[(ii)] Given $X\in \Gee_\xi$\,, $Y\in \Gee_\alpha$ with $\xi<\alpha$, 
	the space $X\vdot\rho_\alpha$ is not almost homeomorphic to $Y$.
\end{enumerate}
We start with $\Gee_0=\Ef_0$. 
Fix $\beta>0$ and suppose $\Gee_\xi$ has been defined for every $\xi<\beta$. 
Let
$$\Aaa_\xi = \setof{X\in\Ef_\beta}
{(\exists\; Y\in\Ef_\xi)\;( X\text{ is almost homeomorphic to }Y\vdot\rho_\beta \, )}.$$
Observe that $|\Aaa_\xi|\loe\kappa_\xi^+$, because $|\Ef_\xi|=\kappa_\xi^+$ 
and, by Lemma~\ref{lemma-3.3},
no two elements of $\Ef_\beta$ are almost homeomorphic. 
Define
$$\Gee_\beta =\Ef_\beta \setminus\bigcup_{\xi<\beta}\Aaa_\xi.$$
Then $|\Gee_\beta|=\kappa_\beta^+$, because 
$|\bigcup_{\xi<\beta}\Aaa_\xi|
\loe\sup_{\xi<\beta}\kappa_\xi^+ 
< \kappa_\beta^+
=|\Ef_\beta|$. 
Hence (i) holds. 
By the definition of $\Gee_\beta$ also (ii) holds.
Finally set $\Gee=\bigcup_{\alpha<\cf\kappa}\Gee_\alpha$. 
Then $|\Gee|=\kappa$, by (i). 
Define
$\Haa = \setof{X\vdot\kappa}{X\in\Gee}$. 
We claim that $\Haa$ satisfies $\star_1(\kappa)$ -- $\star_6(\kappa)$ with 
$\gamma = \kappa$. 
Recall that $\gamma+\kappa = \kappa+\kappa=r(\kappa)$. 
Observe that for $X\in \Gee$, $\rk(X\vdot\kappa) = \kappa$ 
and $e(X\vdot\kappa)$ is the maximal element of $X\vdot\kappa$. 
Also $X\in\calw$ and $|X\vdot\kappa|=\kappa$. 
Thus $\star_i(\kappa)$ holds for $i=1,2,3,4$. 
Condition $\star_6(\kappa)$ follows from the induction hypothesis, 
since every neighborhood of $e(X\vdot\kappa)$ contains a copy of $X$. 
It remains to show $\star_5(\kappa)$. 
For fix $X_0, X_1\in\Gee$ with $X_0\ne X_1$. 
If $X_0,X_1\in\Gee_\alpha\subseteq\Ef_\alpha$ for some $\alpha$ then $X_0$ 
and $X_1$ are not almost homeomorphic by the induction hypothesis, 
therefore so are $X_0\vdot\kappa$ and $X_1\vdot\kappa$, 
by Lemma~\ref{lemma-3.3}. 
Now assume $X_0\in\Gee_\alpha$ and $X_1\in\Gee_\beta$, where $\alpha<\beta$. Observe that $X_0\vdot\kappa$ is isomorphic (as a lattice) to 
$(X_0\vdot\rho_\beta)\vdot\kappa$. 
This is because $X_0\vdot\rho=\fs{P\cdot \rho}$, 
where $P$ is a poset such that $X_0=\fs{P}$; 
clearly $P\cdot \kappa$ is order isomorphic to $(P\cdot\rho)\cdot\kappa$ 
for any $\rho<\kappa$.
By (ii), $X_0\vdot\rho_\beta$ is not almost homeomorphic to $X_1$. 
Recall that $\rk(X_1)=\rho_\beta=\rk(X_0\vdot\rho_\beta)$. 
Hence, by Lemma~\ref{lemma-3.3}, 
$X_0\vdot\kappa\congtop (X_0\vdot\rho_\beta)\vdot\kappa$ 
is not almost homeomorphic to $X_1\vdot\kappa$. 
This shows $\star_5(\kappa)$ and completes the proof.
\end{pf}
\section{Closed subsets of the Tikhonov plank}
It is clear that every uncountable closed subset of the ordinal space 
$X=\omega_1+1$ is homeomorphic (in fact, order isomorphic) to $X$. 
Thus, $X$ has only one non-metrizable compact topological type. 
More generally, given a natural number $n>0$, the space $(\omega_1+1)\times n$
contains only $n$ topological types of non-metrizable compact spaces.
The situation is different in case of the Tikhonov plank 
$\tplank \eqdef (\omega_1+1)\times(\omega+1)$.
Recall that $\tplank$ is a unitary scattered space of rank $\omega_1+1$. 
Moreover, $\tplank$ is a distributive lattice corresponding to the 
free Boolean algebra over $\omega_1 \coprod \omega$.
\begin{theorem}
\label{thm-4.1}
The Tikhonov plank contains $2^{\aleph_1}$ many pairwise 
non-homeo\-mor\-phic closed subsets.
Dually, the Boolean algebra $\F{\omega_1\coprod \omega}$ 
has $2^{\aleph_1}$ many pairwise non-isomorphic quotients.
\end{theorem}
\begin{pf}
Given a closed cofinal set $C\subseteq\omega_1$, define
$$X(C) \eqdef \Bigl((\omega_1+1)\times\sn\omega\Bigr) \cup 
\bigcup_{\xi\in C\cup\sn{\omega_1}}\Bigl(\sn\xi\times(\omega+1)\Bigr).$$
Clearly, $X(C)$ is a closed unitary subspace of $\tplank$, with rank $\omega_1+1$
and $e(X(C))=\pair{\omega_1}{\omega}$. 
We shall construct $2^{\aleph_1}$ many pairwise non-homeomorphic spaces 
of the form $X(C)$.

Let $E(\omega_1) = \setn{ \gamma < \omega_1 }{ \omega^\gamma = \gamma }$ 
be the set of countable $\eps$-ordinals. 
Notice that $E(\omega_1)$ is a closed subset of $\omega_1$.
For every nonempty subset $A$ of~$E(\omega_1)$ we choose 
a function $\fnn{ f_A }{ \omega_1+1 }{ A }$ such that $\inv{f_A}\gamma$ 
contains uncountably many successor ordinals for each $\gamma \in A$.
Given $\beta<  \omega_1$ define
$$\lambda_{A,\beta} \eqdef \sum_{\mu\in A\cap\beta}f_A(\mu),$$
where $\sum$ means the ordinal sum.
Consider the closed and unbounded subset 
$\widehat{A} \eqdef \setn{ \lambda_{A,\beta} }{ \beta <  \omega_1 }$ of $\omega_1 $
and finally define $X_A \eqdef X(\widehat A)$. 
Fix nonempty sets $A,B\subseteq E(\omega_1)$ such that $A\ne B$. 
We shall show that $X_A$ is not homeomorphic to $X_B$.
For suppose $\map g{X_A}{X_B}$ is a homeomorphism. 
Then $g(\pair{\omega_1}\omega) = g(e(X_A)) = e(X_B)= \pair{\omega_1}\omega$. 
Note that $\pair{\omega_1}n$, where $n<\omega$, are the only points of 
rank $\omega_1$ in both spaces $X_A$, $X_B$. 
Thus there exists a bijection $\map\theta{(\omega+1)}{(\omega+1)}$
such that $g(\pair{\omega_1}n) = \pair{\omega_1}{\theta(n)}$ for $n\loe \omega$.

By the continuity of both $g$ and $g^{-1}$, for every $\eta< \omega_1$ we can find
$\eta_{n,0}$, $\eta_{n,1}$ and $\eta_{n,2}$ in the interval $(\eta , \omega_1)$ so that
$$[\eta_{n,0},\omega_1] {\times} \sngltn{\theta(n)}
\subset \img g{[\eta_{n,1},\omega_1] {\times} \sngltn{n}} 
\subset  [\eta_{n,2},\omega_1] {\times} \sngltn{\theta(n)}.$$
We define $\delta(\eta)$ to be the supremum of all ordinals $\eta_{n,i}$, 
where $n\in\omega$ and $i<3$. 
Clearly, $\eta<\delta(\eta)<\omega_1$. 
Define a sequence $\seqn{ \eps_n }{ n\in\omega }$ 
by $\eps_0=0$ and $\eps_{n+1}=\delta(\eps_n)$. 
Let $\eps=\sup_{n\in\omega}\eps_n$. 
Then
$\img g { \, [\varepsilon,\omega_1] {\times} \sngltn{n} \, }
=  [\varepsilon,\omega_1] {\times} \sngltn{\theta(n)}$.
In particular,
\begin{itemize}
\item[\rm(ii)]
$\img g { \bigl(\, [\varepsilon,\omega_1) 
{\cap} \widehat{A}\,\bigr)  {\times} \sngltn{\omega} }
=  \bigl(\, [\varepsilon,\omega_1) {\cap} \widehat{B} \, \bigr) {\times} \sngltn{\omega}$ and
\smallskip
\\
$\img g { \bigl( \, [\varepsilon,\omega_1) {\cap} \widehat{A} \, \bigr) {\times} \omega }
=  \bigl(\, [\varepsilon,\omega_1) {\cap} \widehat{B}\,  \bigr) {\times} \omega$.
\end{itemize}
By (ii) we can define $\map h{[\eps,\omega_1]}{[\eps,\omega_1]}$ 
by setting $g(\pair\mu\omega) = \pair{h(\mu)}\omega$. 
Then $h$ is a homeomorphism of ${[\eps,\omega_1]}$ onto itself.

Now assume $\gamma \in A\setminus B$. 
Recall that $\inv{f_A}\gamma$ contains uncountably many successor ordinals. 
Fix such an ordinal $\alpha$ above $\eps$. 
So $\alpha=\delta+1>\eps$ and $f_A(\alpha)=\gamma$. 
Recall that 
$\lambda_{A,\delta}$\,, $\lambda_{A,\delta+1}$ are members
of $\widehat{A}$ and $\lambda_{A,\alpha+1} = \lambda_{A,\alpha} + \gamma$.
Notice that the interval $[\lambda_{A,\delta}, \lambda_{A,\delta+1}]$ is order isomorphic
(in particular, homeomorphic) to $\gamma+1$. 
Because $\gamma$ is an $\eps$-ordinal,
$$
\rfs{rk}^{X_A}( \pair{\lambda_{A,\delta+1} }{\omega})
= \rfs{rk}^{\omega_1{+}1}(\lambda_{A,\delta+1} ) = \gamma.$$
Since $\lambda_{A,\delta+1}\in\widehat A$, necessarily 
$h(\lambda_{A,\delta+1})\in\widehat B$. 
That is $h(\lambda_{A,\delta+1}) = \lambda_{B,\mu}$ for some $\mu \in \widehat{B}$.
Observe that $[\varepsilon,\omega_1{+}1)$ is a neighborhood of 
both $\lambda_{A,\delta+1}$, $\lambda_{B,\mu}$ 
and that $\lambda_{A,\delta+1}$ is an isolated point in $\widehat{A}$.
Thus the same holds for $\lambda_{B,\mu}$ 
in $\widehat{B} \cap  [\varepsilon,\omega_1{+}1)$.
Hence $\mu$ is a successor, i.e. $\mu=\nu+1$.
But, since $\gamma\notin B$, there is no successor~$\mu$ such that
$\rfs{rk}^{X_B}( \pair{\lambda_{B,\mu} }{\omega}) 
= \rfs{rk}^{\omega_1{+}1}(\lambda_{B,\mu} ) = \gamma$.
This is a contradiction.
\end{pf}

\begin{remark}
A slight modification of the above proof shows that for distinct nonempty 
sets $A,B\subseteq E(\omega_1)$ spaces $X_A$, $X_B$ are not almost homeomorphic. 

Note that each of the spaces $X_A$ contains a copy of $\tplank$, 
therefore it is not retractive. 
Consequently, the algebra of clopen subsets of $X_A$ is not embeddable into any
interval algebra.
\end{remark}
\section{Closed sublattices of $\bm{(\omega_1+1)^2}$}
The purpose of this section is to describe topological types of closed sublattices 
of $(\omega_1+1)^2$. 
In order to simplify some statements, by a {\it sublattice} of a lattice $L$ 
we mean a subset which is closed under meet and join, 
not necessary containing $0_L$ and $1_L$. 
We set 
$$\gplank_\alpha \eqdef (\omega_1+1)\times (\alpha+1)$$
and
$$\gplank \eqdef \gplank_{\omega_1}
= (\omega_1+1)^2\,,\qquad \triangle \eqdef \setof{\pair xy\in\gplank}{x\loe y}\,.$$
We shall prove below that an uncountable closed sublattice of $\gplank$ 
is homeomorphic to one of the above defined spaces. 
The above list of spaces cannot be essentially reduced: 
the only homeomorphic pairs are $\gplank_\alpha$ and $\gplank_\beta$, 
where $\alpha,\beta<\omega_1$ are such that $(\alpha+1)\congtop(\beta+1)$. 
Of course $\gplank_\alpha$ for $\alpha<\omega_1$ is 
neither homeomorphic to $\gplank$ nor to $\triangle$, 
because of the Cantor-Bendixson rank. 
In order to see that $\triangle$ is not homeomorphic to $\gplank$ 
observe that $\cbd{\omega_1}\triangle\congtop \omega_1+1$ 
while $\cbd{\omega_1}{\gplank}\congtop \omega_1+1+\omega_1^*$, 
where $\omega_1^*$ denotes $\omega_1$ with the reversed ordering.

Given topological spaces $X,Y$, we denote by 
$X\oplus Y$ their topological (disjoint) sum. 
Notice that $\gplank_\alpha\oplus\gplank_\beta \congtop \gplank_{\alpha+1+\beta}$ 
and $\gplank_\alpha\oplus \triangle \congtop \triangle$ for a countable ordinal $\alpha$. 
The latter follows from the fact that $\triangle$ contains clopen sets homeomorphic 
to $\gplank_\theta$ for every $\theta<\omega_1$ 
and $\gplank_\alpha\oplus\gplank_\theta\congtop\gplank_\theta$, 
whenever $\theta>\alpha$ is indecomposable. 
Notice also that if $X=\gplank_\alpha$ for some $\alpha\loe\omega_1$ 
or $X=\triangle$ then $(\rho+1)\oplus X\congtop X$ for every countable ordinal $\rho$.
We shall use these observations below.
\begin{lemma}
\label{lemma-5.1}
Let $K$ be a sublattice of $\gplank$ and assume that $\alpha_0\loe\alpha_1$, $\beta_0\loe\beta_1$ are such that $\pair{\alpha_0}{\beta_1}, \pair{\alpha_1}{\beta_0}\in K$. Let $A=\setof{\xi\in[\alpha_0,\alpha_1]}{\pair{\xi}{\beta_0}\in K}$ and $B=\setof{\eta\in[\beta_0,\beta_1]}{\pair{\alpha_0}\eta\in K}$. Then
$$K\cap([\alpha_0,\alpha_1]\times[\beta_0,\beta_1]) = A\times B.$$
\end{lemma}
\begin{pf}
Fix $\pair ab \in A\times B$. 
Then $\alpha_0\loe a$, $\beta_0\loe b$ and $\pair a{\beta_0}, \pair{\alpha_0}b\in K$,
therefore $\pair ab = \pair a{\beta_0} \join \pair{\alpha_0}b \in K$. 
Hence $A\times B\subseteq K$.

Now fix $\pair xy\in K$ such that $\alpha_0\loe x\loe \alpha_1$ 
and $\beta_0\loe y\loe \beta_1$. 
Then $\pair x{\beta_0} = \pair xy \meet \pair{\alpha_1}{\beta_0} \in K$, 
so $x\in A$. Similarly $y\in B$.
\end{pf}
Of course the above lemma is valid for sublattices of products of two arbitrary chains.
\begin{lemma}
\label{lemma-5.2}
Let $K$ be a closed uncountable sublattice of $\gplank_\theta$, 
where $\theta<\omega_1$. 
Then $K$ is homeomorphic to $\gplank_\alpha$ for some $\alpha\loe\theta$.
\end{lemma}
\begin{pf}
Define $\phi(\eta)=\sup\setof{\xi<\omega_1}{\pair\xi\eta\in K}$. 
Let $\rho<\omega_1$ be such that $\phi(\eta)<\rho$ whenever $\phi(\eta)<\omega_1$. 
Let $\beta_0,\beta_1$ be the minimal and the maximal ordinals respectively 
such that $\phi(\beta_0)=\phi(\beta_1)=\omega_1$. 
The maximal one exists by compactness.
Choose $\alpha_0>\rho$ so that $\pair{\alpha_0}{\beta_1}\in K$ and 
let $\alpha_1=\omega_1$. 
Apply Lemma~\ref{lemma-5.1}. 
We conclude that $K_1 = K\cap [\alpha_0,\omega_1]\times[\beta_0,\beta_1]$ 
equals $A\times B$, where the order type of $A$ is necessarily $\omega_1+1$ 
and the order type of $B$ is $\alpha+1$ for some $\alpha\loe\theta$. 
Now the remaining part, that is $K_0 = K\cap ([0,\alpha_0)\times(\theta+1))$ is countable,
therefore homeomorphic to an ordinal. 
By increasing $\alpha_0$, we may assume that $K_0$ is clopen in $K$. 
Hence $K_0\cup K_1$ is homeomorphic to $K_1=A\times B\congtop \gplank_\alpha$.
\end{pf}
\begin{theorem}
\label{thm-5.3}
Let $K$ be an uncountable closed sublattice of $\gplank$. 
Then $K$ is homeomorphic to one of the following spaces: 
$\triangle$, $\gplank$ and $\gplank_\alpha$, where $\alpha<\omega_1$.
\end{theorem}
\begin{pf}
Define $T=\setof{\alpha<\omega_1}{\pair{\alpha}{\omega_1}\in K}$ 
and $R=\setof{\beta<\omega_1}{\pair{\omega_1}{\beta}\in K}$. 
We consider the following cases:
\begin{enumerate}
	\item[(1)] Both sets $T$ and $R$ are nonempty.
	\item[(2)] Both $T$ and $R$ are empty.
	\item[(3)] Exactly one of the sets $T,R$ is nonempty.
\end{enumerate}
We first deal with case (1). Fix $\alpha\in T$, $\beta\in R$. 
By Lemma~\ref{lemma-5.1}, 
$K\cap([\alpha,\omega_1]\times[\beta,\omega_1])=A\times B$ for some sets $A,B$. 
Let $\delta=\max\{\alpha,\beta\}$.
Define
$$U=K\cap([0,\omega_1]\times[0,\delta]),\quad V=K\cap([0,\delta]\times(\delta,\omega_1]),\quad W=K\cap(\delta,\omega_1]^2.$$
Then $U,V,W$ form a partition of $K$ into clopen sublattices. 
By Lemma~\ref{lemma-5.2}, each of the sets $U$, $V$ is either countable 
or homeomorphic to $\gplank_\xi$ for some $\xi<\omega_1$. 
The same applies to $U\cup V=U\oplus V$. 
Further, $W=A'\times B'$, where $A'=A\cap(\delta,\omega_1]$ 
and $B'=B\cap(\delta,\omega_1]$. 
Thus $W$ is either countable or homeomorphic to $\gplank_\eta$ 
for some $\eta\loe\omega_1$. 
We conclude that $K=U\oplus V\oplus W$ is homeomorphic 
to $\gplank_\gamma$ for some $\gamma\loe\omega_1$.
\medskip

We now proceed to case (2). 
We claim that (2) implies $K\congtop\omega_1+1$.
Define $r(\xi)=\sup\setof{\beta}{\pair\xi\beta\in K}$ 
and $u(\eta)=\sup\setof{\alpha}{\pair\alpha\eta\in K}$ for $\xi,\eta<\omega_1$.
Both $r(\xi)$ and $u(\eta)$ are countable ordinals, because $T=R=\emptyset$.
Observe that for every $\delta<\omega_1$ the set $K\cap(\delta,\omega_1]^2$ is 
uncountable; in particular $\pair{\omega_1}{\omega_1}$ is an accumulation point of $K$.
Indeed, $K\setminus (\delta,\omega_1]^2\subseteq(\gamma + 1)^2$, where
$$\gamma = \max \Bigl(\sup_{\xi\loe\delta} r(\xi), \sup_{\eta\loe\delta}u(\eta)\Bigr).$$

Fix a big enough regular cardinal $\chi>\aleph_1$ and fix a continuous chain 
$\seqn{M_\alpha}{\alpha<\omega_1}$ of countable elementary substructures 
of $\pair{H(\chi)}{\in}$ such that $K\in M_0$. 
Let $\theta_\alpha = \omega_1\cap M_\alpha$. 
By elementarity, $\theta_\alpha$ is indecomposable 
and $\ln(\theta_\alpha) = \theta_\alpha$. 
Let $K_\alpha = K\cap [0,\theta_\alpha]^2$. 
Also by elementarity we deduce that $\pair{\theta_\alpha}{\theta_\alpha}\in K_\alpha$.
Recall that $\rk^K(x)$ denotes the Cantor-Bendixson rank of $x$ in $K$. 
By elementarity $\rk^{K_\alpha}(x) \in M_\alpha$ whenever $x \in K_\alpha$.
Observe that 
$K_\alpha \setminus \sn{ \pair{\theta_\alpha}{\theta_\alpha} } \subseteq M_\alpha$.
This is because, again by elemetarity, $r(\xi)<\theta_\alpha$ and $u(\eta)<\theta_\alpha$
for every $\xi,\eta<\theta_\alpha$. 
It follows that $\rk(K_\alpha)\loe\theta_\alpha$. 
On the other hand, for every $\xi<\theta_\alpha$ there is an element $x$ of $K$ 
with $\rk^K(x)\goe\xi$, therefore by elementarity there is also such an element in 
$M_\alpha$. 
That is, for every $\xi<\theta_\alpha$ there is an element $x$ of $K_\alpha$ 
such that $\rk^{K_\alpha}(x) \in M_\alpha$ and $\rk^{K_\alpha}(x)\goe\xi$.
This shows that $\rk(K_\alpha)=\theta_\alpha$ and 
$e(K_\alpha)=\pair{\theta_\alpha}{\theta_\alpha}$. 
Thus $K_\alpha\congtop \theta_\alpha+1$.

It is now straight how to find a homeomorphism $\map hK{\omega_1+1}$. 
We construct inductively a sequence of homeomorphisms 
$\map {h_\alpha}{K_\alpha}{\theta_\alpha+1}$ such that 
$h_\beta$ extends $h_\alpha$ whenever $\beta>\alpha$.
At a successor stage, notice that $K_\alpha$ is clopen in $K_{\alpha+1}$ 
and $h_\alpha$ can be extended to a homeomorphism 
$\map{h_{\alpha+1}}{K_{\alpha+1}}{\theta_{\alpha+1}+1}$, 
because $\theta_{\alpha+1}$ is indecomposable.
Suppose now that $\delta$ is a limit ordinal and $h_\alpha$ has been defined for every 
$\alpha<\delta$. 
Observe that 
$K_\delta = \sn{\pair{\theta_\delta}{\theta_\delta}}\cup\bigcup_{\alpha<\delta}K_\alpha$, 
by the continuity of the chain $\sett{M_\alpha}{\alpha<\omega_1}$ and by the above remarks.
Let $h_\delta$ be the extension of $\bigcup_{\alpha<\delta}h_\alpha$ defined by 
$h_\delta(\pair{\theta_\delta}{\theta_\delta})=\theta_\delta+1$. 
We only need to verify the continuity of $h_\delta$ at $\pair{\theta_\delta}{\theta_\delta}$.
Fix a neighborhood $U=(\rho,\theta_\delta]$ of $\theta_\delta$ in $\theta_\delta+1$. 
We may assume that $\rho = \theta_\alpha$ for some $\alpha<\delta$.
Recall that $\img{h_\delta} {K_\alpha} = \theta_\alpha+1$,
so $\inv{h_\delta}U=K_\delta\setminus K_\alpha$, which is clopen in $K_\delta$. 
This shows the continuity of $h_\delta$. 
Finally $h \eqdef h_{\omega_1}$ is the desired homeomorphism of $K$ 
onto $\omega_1+1$.
\medskip

We are left with case (3). 
We assume that $T=\emptyset$ and $R\ne\emptyset$ (both possibilities are symmetric). 
Define the function $r$ as before, i.e. $r(\xi)=\sup\setof{\beta}{\pair{\xi}{\beta}\in K}$ 
for $\xi<\omega_1$. 
Clearly $r(\xi)<\omega_1$, because $T=\emptyset$. 
Suppose that $\delta = \sup_{\xi<\omega_1}r(\xi)$ is countable. 
Let $U=K\cap [0,\omega_1]\times[0,\delta]$ and 
$V = K\cap \bigl( \sn{\omega_1}\times(\delta,\omega_1] \bigr)$. 
Clearly, $U$, $V$ are disjoint clopen sublattices of $K=U\cup V$. 
Again by Lemma~\ref{lemma-5.2}, $U$ is either countable or homeomorphic 
to $\gplank_\rho$ for some $\rho\loe\delta$. 
Finally, $V$ is either countable or homeomorphic to $\omega_1+1$, 
so we conclude that $K\congtop \gplank_\gamma$ for some $\gamma\loe\delta+1$.

Assume now that the function $r$ is unbounded,
i.e. $\sup_{\xi<\omega_1}r(\xi) = \omega_1$. 
We claim that $K\congtop\triangle$. 
First we show that the order-type of $R$ is $\omega_1$.
Indeed let $\nu \in R$. 
Choose $\xi<\omega_1$ such that $r(\xi) > \nu$
and $\nu' \in [\nu , \omega_1)$ such that $\pair{\xi}{\nu'} \in R$.
Then $\pair{\omega_1}{\nu'} = \pair{\xi}{\nu'} \join \pair{\omega_1}{\nu} \in R$.
That is $\nu' \in R$.
Hence $R$ is unbounded in $\omega_1$.

Next, as in case (2), we fix a continuous chain $\seqn{M_\alpha}{\alpha<\omega_1}$ 
of countable elementary substructures of a big enough $\pair{H(\chi)}{\in}$, 
such that $K\in M_0$. 
Let $\theta_\alpha = \omega_1\cap M_\alpha$ and 
let $K_\alpha=K\cap[0,\theta_\alpha]^2$. 
Further, let $\triangle_\alpha = \triangle\cap[0,\theta_\alpha]^2$ and 
let $R_\alpha = R\cap\theta_\alpha$.
Since the order-type $R$ is $\omega_1$,
by elementarity, the order type of $R_\alpha$
is $\theta_\alpha$.

Fix $\alpha<\omega_1$ and fix $\eta\in R_\alpha$. 
We claim that for every $\xi<\omega_1$ there is $\xi'\in[\xi,\omega_1)$ 
such that $\pair{\xi'}\eta\in K$. 
Indeed, since $r$ is unbounded, we may find an arbitrary large countable ordinal $\xi'$
such that $r(\xi')>\eta$. 
Let $\pair{\xi'}\rho\in K$ be such that $\rho\goe\eta$. 
Then $\pair{\xi'}\eta = \pair{\xi'}\rho\meet\pair{\omega_1}\eta \in K$. 
By elementarity, for every $\xi<\theta_\alpha$ there is $\xi'\in[\xi,\theta_\alpha)$ 
such that $\pair{\xi'}\eta\in K_\alpha$. 
This shows that $\pair{\theta_\alpha}\eta$ is an accumulation point 
of $K_\alpha\cap([0,\theta_\alpha)\times\sn\eta)$. 
For every $\xi<\theta_\alpha$ there is an element $x$ of 
$K\cap([0,\omega_1) \times \sn\eta)$
with $\rk^K(x)\goe\xi$, therefore by elementarity there is also such an element in 
$M_\alpha$. 
Since $\theta_\alpha$ is indecomposable and $\rfs{ln}(\theta_\alpha) = \theta_\alpha$,
we conclude that 
$\rk^{K_\alpha}(\pair{\theta_\alpha}{\eta})\goe\theta_\alpha$. 
It follows that
$$\cbd{\theta_\alpha}{K_\alpha} 
= K_\alpha \cap \bigl( \sn{\theta_\alpha}\times\omega_1 \bigr)
= \sn{\theta_\alpha}\times \bigl( R_\alpha\cup\sn{\theta_\alpha} \bigr) \, .$$
A similar fact applies to $\triangle_\alpha$, 
that is $\cbd{\theta_\alpha}{\triangle_\alpha}=\sn{\theta_\alpha}\times(\theta_\alpha+1)$. 
Recall that the order type of $R_\alpha$ is $\theta_\alpha$. 
It follows that $K_\alpha$ is homeomorphic to $\triangle_\alpha$ and every 
homeomorphism between these spaces maps $\sn{\theta_\alpha}\times R_\alpha$ 
onto $\sn{\theta_\alpha}\times\theta_\alpha$. 
Given a homeomorphism $\map f{K_\alpha}{\triangle_\alpha}$, we shall define 
$\map{f^*}{R_\alpha}{\theta_\alpha}$ to be the unique map satisfying the equation 
$f(\pair{\theta_\alpha}{\eta}) = \pair{\theta_\alpha}{f^*(\eta)}$.

We construct a sequence of homeomorphisms $\map{h_\alpha}{K_\alpha}{\triangle_\alpha}$ satisfying the following conditions:
\begin{enumerate}
	\item[(i)] $h_\beta$ extends $h_\alpha$ whenever $\beta>\alpha$.
	\item[(ii)] If $\alpha<\beta$ then $f_\beta(\pair\xi\eta) = \pair\xi{f_\alpha^*(\eta)}$ for every $\pair\xi\eta\in K_\beta\setminus K_\alpha$ such that $\eta<\theta_\alpha$.
\end{enumerate}
Condition (ii) requires some explanation. 
Fix $\alpha<\beta<\omega_1$. 
Define
$$K^\beta_\alpha 
= K_\beta \cap \bigl( [\theta_\alpha,\theta_\beta]\times(\theta_\alpha,\theta_\beta] \bigr)
\text{ \ and \ }
B^\beta_\alpha = K_\beta\setminus \bigl( K_\alpha\cup K^\beta_\alpha \bigr) \, .$$
Similarly, define
$$\triangle^\beta_\alpha 
= \triangle_\beta \cap \bigl( [\theta_\alpha,\theta_\beta]\times(\theta_\alpha,\theta_\beta] \bigr) 
\text{ \ and \ }
\square^\beta_\alpha 
= \triangle_\beta\setminus \bigl( \triangle_\alpha\cup\triangle^\beta_\alpha \bigr) \,.$$
Clearly, $\square^\beta_\alpha$ is a rectangle of the form 
$(\theta_\alpha,\theta_\beta] \times [0,\theta_\alpha]$. 
We claim that $B^\beta_\alpha$ is also a rectangle, whose horizontal side 
has order type $\theta_\beta$ and whose vertical side is $R_\alpha\cup\sn{\theta_\alpha}$
(so its order type is $\theta_\beta{+}1$ as well).
Let $\pi=\min R$. 
Then $\pi\in M_0\subseteq M_\beta$, so again by elementarity we deduce that 
$\pair{\theta_\beta}\pi$ is an accumulation point of 
$K_\beta\cap (\theta_\beta\times\sn\pi)$. 
In particular, $\pair{\theta_\beta}\pi\in K_\beta$. 
Also $\pair{\theta_\alpha}{\theta_\alpha}\in K_\beta$, so by Lemma~\ref{lemma-5.1},
$$K_\beta\cap([\theta_\alpha,\theta_\beta]\times[\pi,\theta_\beta]) = C\times D,$$
where $D \eqdef R_\alpha\cup\sn{\theta_\alpha}$ and the order type of $C$ is necessarily 
$\theta_\beta+1$. 
Thus $B^\beta_\alpha = \bigl(C\setminus\sn{\theta_\alpha} \bigr) \times D$ is naturally
homeomorphic to $\square^\beta_\alpha$ and condition (ii) makes sense.

We start with any homeomorphism $\map{f_0}{K_0}{\triangle_0}$, which exists 
because both spaces are unitary of rank $\theta_0+\theta_0$. 
Fix $\beta>0$ and suppose that homeomorphisms $f_\alpha$ have been defined 
for all $\alpha<\beta$. 
Suppose first that $\beta=\alpha+1$. 
By the above remarks, there is a unique way to define $f_{\alpha+1}$ on 
$B^{\alpha+1}_\alpha$ so that (ii) holds. 
Also, $K^{\alpha+1}_\alpha$ is homeomorphic to $\triangle^{\alpha+1}_\alpha$ 
so it is possible to extend $f_\alpha$ to $f_{\alpha+1}$ and any such extension 
maps $\pair{\theta_\beta}{\theta_\beta}$ onto $\pair{\theta_\beta}{\theta_\beta}$.

Suppose now that $\beta$ is a limit ordinal. 
The set $D=\bigcup_{\xi<\beta}K_\xi$ is dense in $K_\beta$ and 
$K_\beta\setminus D 
= \bigl( R_\beta\cup\sn{\theta_\beta}\bigr)\times\sn{\theta_\beta \bigr} $. 
We need to define $f_\beta\supseteq \bigcup_{\alpha<\beta}f_\alpha$. 
Fix $\eta<\theta_\beta$. 
Then $\eta<\theta_\alpha$ for some $\alpha<\beta$. 
Let $f_\beta(\pair{\theta_\beta}{\eta}) = \pair{\theta_\beta}{f^*_\alpha(\eta)}$. 
Condition (ii) says that $f^*_\xi(\eta)=f^*_\alpha(\eta)$ for every $\xi\in[\alpha,\beta)$, 
thus (ii) holds also for $\beta$. 
Since $f^*_\alpha$ is continuous, we deduce that $f_\beta$ is continuous 
at $\pair{\theta_\beta}{\eta}$. 
Finally set $f_\beta(\pair{\theta_\beta}{\theta_\beta}) = \pair{\theta_\beta}{\theta_\beta}$. 
It is clear that $f_\beta$ is continuous at $\pair{\theta_\beta}{\theta_\beta}$.

Finally, the map $\map{f_{\omega_1}}K\triangle$ which is the unique extension of 
$\bigcup_{\alpha<\omega_1}f_\alpha$, is a homeomorphism. 
Its continuity follows by the same argument as in the limit stage of the construction. 
This completes the proof.
\end{pf}
\begin{corollary}
Let $P=P_0\cup P_1$ be an uncountable poset such that both $P_0$, $P_1$ are 
well ordered chains of order type $\loe\omega_1$. 
Then the Boolean algebra $\F P$ is isomorphic to one of the following algebras:
$\F{\omega_1}$, $\F{\omega_1\uplus\omega_1}$, $\F{\omega_1\uplus \alpha}$, 
where $\alpha$ is a countable ordinal, and $\F{\omega_1\times2}$, 
where $\omega_1\times2$ is endowed with the coordinatewise order.
\end{corollary}
\begin{pf}
Since $P$ is uncountable, we may assume that the order type of $P_0$ is $\omega_1$. 
Let $\alpha$ be the order type of $P_1$ and let $\map f{\omega_1\uplus\alpha}P$ be 
the natural order preserving surjection.
By duality, $f$ corresponds to a topological lattice embedding of $\fs P$ into 
$\fs{\omega_1\uplus\alpha}$. 
Now the corollary follows from Theorem~\ref{thm-5.3}, 
noticing that $\triangle$ is isomorphic to $\fs{\omega_1\times2}$.
\end{pf}

\end{document}